\RequirePackage{ifpdf}
\ifpdf % We are running pdfTeX in pdf mode
\documentclass[pdftex]{sigma}
\else
\documentclass{sigma}
\fi

\def\1{{\bf 1}}
\def\2{{\bf 2}}
\def\3{{\bf 3}}
\def\j{{\bf j}}
\def\k{{\bf k}}
\def\m{{\bf m}}
\def\C{{\cal C}}
\def\D{{\cal D}}
\def\E{{\cal E}}
\def\F{{\cal F}}
\def\p{\partial}
\def\P{\mathbb{P}}
\def\R{\mathbb{R}}
\def\W{{\cal W}}

\begin{document}

\allowdisplaybreaks

\renewcommand{\thefootnote}{$\star$}

\renewcommand{\PaperNumber}{102}

\FirstPageHeading

\ShortArticleName{Singularity Classes of Special 2-Flags}

\ArticleName{Singularity Classes of Special 2-Flags\footnote{This paper is a
contribution to the Special Issue ``\'Elie Cartan and Dif\/ferential Geometry''. The
full collection is available at
\textit{}\href{http://www.emis.de/journals/SIGMA/Cartan.html}{http://www.emis.de/journals/SIGMA/Cartan.html}}}

\Author{Piotr MORMUL}

\AuthorNameForHeading{P. Mormul}

\Address{Institute of Mathematics, Warsaw University,
2 Banach Str., 02-097 Warsaw, Poland}

\Email{\href{mailto:mormul@mimuw.edu.pl}{mormul@mimuw.edu.pl}}

%\URLaddress{\url{http://www.mimuw.edu.pl/~mormul/}}

\ArticleDates{Received April 16, 2009, in f\/inal form October 30, 2009;  Published online November 13, 2009}

\Abstract{In the paper we discuss certain classes of vector distributions
in the tangent bundles to manifolds, obtained by series of applications
of the so-called generalized Cartan prolongations (gCp). The classical
Cartan prolongations deal with rank-2 distributions and are responsible
for the appearance of the Goursat distributions. Similarly, the so-called
special multi-f\/lags are generated in the result of successive applications
of gCp's. Singularities of such distributions turn out to be very rich,
although without functional moduli of the local classif\/ication.
The paper focuses on special 2-f\/lags, obtained by sequences
of gCp's applied to rank-3 distributions.
A stratif\/ication of germs of special 2-f\/lags of all lengths into
{\it singularity classes} is constructed. This stratif\/ication
provides invariant geometric signif\/icance to the vast family
of local polynomial pseudo-normal forms for special 2-f\/lags
introduced earlier in [Mormul P., {\it Banach Center Publ.}, Vol.~65, Polish Acad. Sci., Warsaw, 2004, 157--178]. This is the main
contribution of the present paper. The singularity classes
endow those multi-parameter normal forms, which were obtained
just as a by-product of sequences of gCp's, with a geometrical
meaning.}

\Keywords{generalized Cartan prolongation; special multi-f\/lag;
special 2-f\/lag; singularity class}

\Classification{58A15; 58A17; 58A30}

\renewcommand{\thefootnote}{\arabic{footnote}}
\setcounter{footnote}{0}

\section{Introduction and main theorem}\label{IMT}

The aim of the current paper is to present a new and rather rich
stratif\/ication of singularities of (special) 2-f\/lags which naturally
generalize 1-f\/lags. Before doing that, it will be useful to brief\/ly
recall 1-f\/lags and their singularities. These are, in the contemporary
terminology, rank-2 and corank $\ge 2$ subbundles $D \subset TM$ in
the tangent bundle to a smooth manifold $M$, together with the tower
of consecutive {\it Lie squares} $D \subset [D,D] \subset
[[D,D],[D,D]] \subset \cdots$ satisfying the property
that the linear dimensions of tower's members are $2,3,4,\dots$
at {\it every} point in~$M$. (In $(\dim M - 2)$ steps the tower
reaches the full tangent bundle $TM$.) These objects had emerged in
the papers \cite{Engel,vonWeber,Cartan1914} and were later popularized
in a book by Goursat in the 1920s. In the result, such distributions
$D$ are now called the {\it Goursat distributions}, or sometimes the
{\it Cartan--Goursat distributions}.
The respective f\/lags are called the {\it Goursat flags}. Although this
def\/inition is quite restrictive, still such f\/lags exist in all lengths.
Indeed, for every $s \ge 2$, the canonical contact system $\C^s$
(the jet bundle or the Cartan distribution in the terminology of
\cite{Krasilshchik}) on the jet space $J^s(1,1)$ is a Goursat
distribution of corank $s$; its f\/lag has length $s$. However, each
distribution $\C^s$ is homogeneous because its germs at every two
points are equivalent by a local dif\/feomorphism of $J^s(1,1)$.
Therefore, these contact systems have no singularities.
It should be noted that nowadays the contact systems
on $J^s(1,1)$ are {\it also} known under the name
`Goursat normal forms' and are characterized as such
in \cite{Bryant1991} (Theorem~5.3 in Chapter~II).

For a very long time it had not been known whether Goursat f\/lags
locally featured any other geometry than that of the systems $\C^s$.
An af\/f\/irmative answer was given only in 1978 by Giaro, Kumpera
and Ruiz in dimension~5, and slightly later, in~\cite{Kumpera1982},
in all dimensions $\ge 5$. Much later a geometric systematization of
those f\/indings appeared in \cite{Montgomery2001}. Namely, Montgomery
and Zhitomirskii def\/ined {\it Kumpera--Ruiz classes} of germs of
Goursat f\/lags (KR classes for short) in every f\/ixed length
$s \ge 2$. The number of them in length $s$ is $2^{s-2}$.
In fact, it is natural to encode those classes by words of
length $s$ over the alphabet $\{1,2\}$.  The words start with
two 1's. In the $i$-th place, where $3 \le i \le s$, one writes
1 or 2 depending on whether the condition (GEN) from p.~466
in~\cite{Montgomery2001} holds true for that $i$ or not.
This specif\/ication of the way in which one puts the numbers is
purely geometrical and means that either certain two (invariantly
def\/ined) lines in a~plane, related with the corank-$i$ member of
the f\/lag at the reference point, are dif\/ferent or merge into one
line. Nearly immediately those classes appeared to perfectly match
the $2^{s-2}$ branches in the tree of Kumpera--Ruiz [pseudo]normal
forms for germs of Goursat distributions of corank~$s$ constructed
in~\cite{Kumpera1982}. (Those were couples of {\it polynomial}
vector f\/ields with only f\/inite number of {\it real parameters}.
The construction of those f\/ields had much in common with a KR class
to which the relevant germ belonged.) It was a departure point for
an entirely new, full-scale theory of Goursat f\/lags developed in
recent years.

Let us emphasize the key fact which has motivated the present article.
The following two seemingly distant aspects of the theory are closely
related:
\begin{enumerate}\itemsep=0pt
\item[--] the local realizations, or KR normal forms constructed
in \cite{Kumpera1982},  and

\item[--] the genuine KR classes of singularities def\/ined
in \cite{Montgomery2001}.
\end{enumerate}
 (The former preceded the latter by 18 years!)

Our objective is to establish an analogous, but going further,
relationship for very natural generalizations of 1-f\/lags,
the so-called {\it special} 2-{\it flags.}

So, to begin with, what are special multi-f\/lags? In the def\/inition
we will use the notion of the {\it Cauchy-characteristic module} (or,
strictly speaking, sheaf of modules) of a distribution $D$, written
$L(D)$ (the Japanese school adheres to the symbol ${\rm Ch}(D)$).
 It consists of all vector f\/ields~$v$ (in the considered category of
smoothness) taking values in $D$ and preserving $D$: $[v,D] \subset D$.

\begin{definition}[special $\boldsymbol{m}$-f\/lags] \label{definition1} We f\/ix a natural
number $m \ge 2$ (called `width').  A rank-$(m+1)$ distribution $D$
on a manifold $M$ generates a special $m$-f\/lag of length $r \ge 1$
on $M$ when
\begin{enumerate}\itemsep=0pt
\item[$\star$] the tower of consecutive Lie squares of $D$
\[
D = D^r \subset D^{r-1} \subset D^{r-2} \subset \cdots \subset D^1
\subset D^0 = TM,
\]
$[D^j,D^j] = D^{j-1}$ for $j = r,r-1,\dots,2,1$, consists
of distributions of ranks, starting from the smallest object $D^r$:
$m+1, 2m+1,\dots, rm + 1, (r+1)m + 1 = \dim M$,

\item[$\star\star$] for $j = 1,2,\dots,r-1$ the Cauchy-characteristic
module $L(D^j)$ of $D^j$ sits already in the smaller object $D^{j+1}$,
$L(D^j) \subset D^{j+1}$, and is regular of corank 1 in $D^{j+1}$
(i.e., such a module of vector f\/ields has its linear dimension
${\rm rk}\, D^{j+1} - 1$ at every point),  while $L(D^r) = 0$,

\item[$\star\star\star$] the biggest f\/lag's member $D^1$ possesses
a corank-1 involutive (i.e., completely integrable) subdistribution,
which we call $F$.
\end{enumerate}
\end{definition}

This def\/inition is slightly more specif\/ic than the original
def\/inition from~\cite{Mormul2004a}. It is, however, equivalent, singling
out precisely the same objects. It emphasizes the Cauchy-characteristic
subf\/lag; compare also the def\/inition of `generalized contact systems
for curves' in \cite{Pasillas}. On the other hand, Definition~\ref{definition1} is,
as it stands, redundant, for condition $\star\star$ is implied by
$\star$ and $\star\star\star$, see Proposition~1.3 in \cite{Adachi},
or Corollary~6.3 in~\cite{Shibuya}. Thus the meaning of `special'
resides in the existence of an involutive corank-1 subdistribution
$F \subset D^1$.\footnote{Such $F$, when exists, is unique and
has more than one geometrical interpretation (see Corollary~\ref{unique}
and Remarks~\ref{remark1} and~\ref{remark2}(b) later on, and also~\cite[p.~165]{Mormul2004a}).
It is worth noting that when $m = 1$ such subdistributions~$F$ also exist
but are {\it not} unique and have no geometrical meaning whatsoever
(this concerns Goursat f\/lags which are not considered in the present
paper).} The involutiveness of $F$ is critical; see Remark~\ref{remark2}(a) in
this respect.

Additional comment to Def\/inition~\ref{definition1}: condition $\star$ alone def\/ines
{\it general} $m$-f\/lags, whose possible geometries are extremely
rich, including functional moduli in the local classif\/ication, see for
instance~\cite{Cartan1910,Agrachev,Zelenko}. It is neatly outbalanced
by conditions $\star\star$ and $\star\star\star$ (in fact, reiterating,
the latter eventually covers the former).

 As for $m = 1$, that time condition $\star\star$ is implied
by just condition $\star$ and so 1-f\/lags appear to be automatically special!
(This is outside the scope of special multi-f\/lags.)

There exist ef\/fective realization techniques producing distributions which
generate special multi-f\/lags of arbitrary width and length. For general $m \ge 2$
they have been constructed in Section~3.3 of~\cite{Mormul2004a} with an essential
use of the so-called generalized Cartan prolongations; see in this respect
Theorem~\ref{not} later on. (Specif\/ically for $m = 2$, for reader's convenience,
these operations are re-def\/ined in Sections~\ref{san} and~\ref{EKRs} of
the present paper.) In the outcome one gets a~vast family of polynomial
[pseudo-]normal forms with many numerical parameters, the so-called Extended
Kumpera--Ruiz normal forms (EKR for short; see Section~\ref{EKRs} for the
explanation of the origin of this name).
Within a given EKR the realizations dif\/fer only by the values of
the numeric parameters that enter that EKR. The classes EKR are encoded
(or: labelled) by words $\j_1\!\dots\j_{r-1}.\j_r$ over $\{\1,\2,\dots,\m,\m+\1\}$ subject to an important limitation called
{\it the least upward jumps rule}. Namely, admissible words should start
with \1 and always a {\it new} but not yet used number should only
{\it minimally} exceed the maximum of previously used numbers:
for $l = 1,2,\dots,r-1$, if $\j_{l+1} > \max (\j_1,
\dots,\j_l)$ then $\j_{l+1} = 1 + \max(\j_1,\dots,\j_l)$.

 For $m = 2$ this rule says that, after starting from~\1, the f\/irst
use of \3 (if any) should occur after the f\/irst use of~\2 (if any). That is
to say, the number~\3 does not appear without number~\2 before it. It is
straightforward to see (Proposition~\ref{how}) that the number of EKR
classes is equal to $\frac{1}{2}\bigl(1 + 3^{r-1}\bigr)$ in every
length $r \ge 1$.

 Instead of the Kumpera--Ruiz normal forms for Goursat f\/lags,
we now have EKR's to ef\/f\/i\-cient\-ly handle special multi-f\/lags. Indeed,
neat polynomial local realizations have been proposed in both settings
$m = 1$ and $m \ge 2$. In the former case it is known (and already mentioned
above) that the KR normal forms faithfully correspond to the KR classes of
singularities put forward in~\cite{Montgomery2001}.  Do, therefore, the
EKR's in the latter case correspond to some partition or stratif\/ication
of the space of all germs of special $m$-f\/lags\footnote{This
was B.~Kruglikov's question asked in 2002.}? Or, in the least, what
could be said specif\/ically in width~2?

Our objective in the present paper is to answer the main question above
af\/f\/irmatively in width $m = 2$ for the rank-3 distributions generating
special 2-f\/lags of arbitrary length.

 Firstly in Section~\ref{san} we construct an analogue of the
KR classes of Goursat f\/lags for special 2-f\/lags, adapting the method of~\cite{Montgomery2001}.
We call the obtained intermediate aggregates of germs of special 2-f\/lags
`sandwich classes', because they directly emanate from the sandwich diagram
for multi-f\/lags (see Section~\ref{SD})\footnote{In this step the
specif\/ication $m = 2$ is not important.}. The sandwich classes are encoded
by such words over $\{1,\underline{2}\}$ which start with 1 and are of
length equal to f\/lag's length. If the length is $r$, then the number of
sandwich classes is $2^{r-1}$ (note the dif\/ference in exponent with the
Goursat case, which is due to the presence of the distribution $F$ in
the sandwich diagram for multi-f\/lags).

 Secondly, we present the key part of the paper in Section~\ref{refin}.
Namely, only in width 2, we ref\/ine the notion of the sandwich class to a
`singularity class'. In fact, a germ $D$ sitting in a~given sandwich class
$\cal S$ of length $r$ is being analyzed both geometrically and Lie-algebraically.
The purpose is to specify all but the f\/irst~\underline{2} (from the left)
in the label of $\cal S$ to~2 or~3, each one independently of the others.
It is the local geometry of the f\/lag of $D$ that decides that choice.
(As for the f\/irst \underline{2}, it is invariably specif\/ied to~2.)
In the outcome a word $j_1.j_2\dots j_r$ over $\{1,2,3\}$, denoted
by~$\W(D)$, is being associated to $D$. Given that the `sandwich' words
start with 1 and the f\/irst~\underline{2} in them is later specif\/ied to~2,
it is clear that $\W(D)$ also satisf\/ies the least upward jumps rule,
exactly as the labels of the EKR classes have done. We mean that $j_1 = 1$
and if $j_{l+1} > \max(j_1,\dots,j_l)$ then $j_{l+1} =
1 + \max(j_1,\dots,j_l)$ for $l = 1,2,\dots,r-1$.

Now the germs having one and the same word $\W(\cdot)$ build up
a given {\it singularity class}. It follows that the partition of all germs
into singularity classes is a ref\/inement of the partition into sandwich
classes, and that the cardinality of that new f\/iner partition is
the same as the cardinality of the EKRs in that length. That is,
$\frac{1}{2}\bigl(1 + 3^{r-1}\bigr)$ in every length $r$
(Proposition~\ref{how}).

 Does one know that all singularity classes are nonempty? More
generally, is there a relationship among the singularity classes and the
classes EKR of concrete realizations of special 2-f\/lags in any given length~$r$? Are singularity classes visible on the level of local polynomial pseudo-normal
forms EKR? It turns out that the answer is `yes' and the EKRs do not forget about
the underlying local f\/lag's geometry concretized by (or: discretized in)
the singularity class. Namely, there holds

\begin{theorem}[main theorem]\label{2003}
Let $D$ be any germ of a rank-$3$ distribution generating a special $2$-flag
of length $r \ge 1$, belonging to a fixed singularity class $j_1.j_2\dots j_r$
$(=\W(D))$. Then the  EKR pseudo-normal forms of $D$ are  uniquely
of the type $\j_1.\j_2\dots \j_r$.
\end{theorem}

Although elusive on the def\/inition level, f\/lag's local invariant -- singularity
class -- acquires a~concrete illustration in this theorem.

\begin{corollary}\label{lass}
The singularity class of a germ of a special $2$-flag which is already given
in an EKR form $\j_1.\j_2\dots \j_r$ is $j_1.j_2\dots j_r$.
That is, slightly abusing notation, $\W(\j_1.\j_2\dots \j_r) = j_1.j_2\dots j_r$.
\end{corollary}

Therefore, Theorem~\ref{2003} {\it additionally} shows that all singularity
classes are non-empty. Whenever one f\/inds an EKR for a germ of special 2-f\/lag,
one inevitably stumbles upon its singularity class. An illustrative example
of retrieving the singularity classes from EKRs is given later in
Appendix~\ref{concr}. Theorem~\ref{2003} is proved in
Section~\ref{pf}.

\section{Generalized Cartan prolongations produce special multi-f\/lags}
\label{ssecond}

In dif\/ferential geometry there exists an important operation, def\/ined in
the papers of \'E.~Cartan and used by him in various situations. It takes
rank-2 vector distributions in the tangent bundles to manifolds, processes
them and yields more complicated rank-2 distributions, living on bigger
manifolds, in the outcome. Nowadays it is called {\it Cartan prolongation}
and can be applied to an arbitrary rank-2 distribution. In the modern
language of \cite[p.~454]{Bryant1993} its def\/inition goes
as follows.

`If $\D$ is a rank-2 distribution on a manifold $M$, then, regarding $\D$
as a vector bundle, we can certainly def\/ine its projectivization $\pi: \P\D
\longrightarrow M$, which is a bundle over $M$ whose typical f\/iber~$\P\D_p$
is the space of 1-dimensional linear subspaces of the 2-dimensional vector
space $\D_p$. Thus, the f\/ibers of $\P\D$ are isomorphic to $\P\R^1$ as
projective 1-manifolds. There is a canonical rank-2 distribution $\D^{(1)}$
on $\P\D$ def\/ined by setting $\D^{(1)}_\xi = (\pi')^{-1}(\xi)$ for each
linear subspace $\xi \subset \D_p$. The distribution $\D^{(1)}$ is called
the {\it $($first$)$ prolongation} of $\D$.'

Its importance stems from a key local {\it structural} theorem, presented
(or, in authors' view, only recalled) in~\cite{Bryant1993}. This theorem deals
with rank-2 distributions that have mild properties of growing (in a Lie algebra
sense which is explained below) neither too slowly nor too quickly. In fact,
certain rank-2 distributions of corank, say, $s$ locally turn out to be,
up to the equivalence of the base manifolds, nothing but the Cartan prolongations
of rank-2 distributions of corank~$s-1$. This shows that such distributions
are constructed simpler than one could expect, and that they have some handy
structure.

 Below, $\D_1$ means the Lie square $[\D,\D]$ of a distribution $\D$,
and $\D_2$ -- the Lie square of~$\D_1$. The foliation~$\F$ is a classical
object closely related to the hypothesis on the def\/icient rank of~$\D_2$
(everywhere~4 instead of~5). In fact, under the hypotheses in Theorem~\ref{cbh},
the Cauchy-characteristic module $L(\D_1)$ (for the def\/inition of $L(\cdot)$,
see the paragraph before Def\/inition~\ref{definition1} earlier on) is a rank-1 subdistribution
of $\D_1$~-- a f\/ield of lines, and $\F$ is the 1-dimensional foliation
tangent to $L(\D_1)$.

\begin{theorem}[Cartan--Bryant--Hsu]\label{cbh}
Let $\D$ be a rank-$2$ distribution on a manifold $M^{s+2}$ and suppose
that $\D_1$ and $\D_2$ have ranks~$3$ and~$4$ respectively. Furthermore,
suppose that there is a~submersion $f: M \rightarrow N^{s+1}$ with the
property that the fibers of  $f$ are the leaves of the canonical foliation~$\F$. Then there exists a unique rank-$2$ distribution $\D'$ on N with the
property that $\D_1 = f^*(\D')$  and, moreover, there exists a canonical
smooth map $f^{(1)}:M \longrightarrow \P\D'$ which is a~local
diffeomorphism, which satisfies $f = \pi\circ f^{(1)}$, and
which satisfies $f^{(1)}_{\ \,*}\D = (\D')^{(1)}$.
\end{theorem}

The ultimate consequence of this impressive theorem is a clear local
construction of Goursat distributions. As simply as it can only be,
Cartan prolongations applied in longer and longer successions produce
(locally) all longer and longer Goursat f\/lags!  Montgomery and Zhitomirskii
summarize the resulting situation in \cite[p.~479]{Montgomery2001}
as follows:
`Every corank $s$ Goursat germ can be found, up to a dif\/feomorphism,
within the $s$-fold prolongation of the tangent bundle to a~surface.
We have called this $s$-fold prolongation the ``monster manifold''.
It is a very tame monster in many respects.'

A recent big contribution \cite{Montgomery2008} of the same authors
demonstrates how eventually fruitful this Cartan-inspired visualisation
of Goursat distributions is.

Returning to [special] multi-f\/lags, an instance of vagueness shrouding
them 10 years ago is the following. The f\/irst of the authors of
\cite{Kumpera2002} wrote, in a personal communication,
in spring of 1999:

\dots {\it multi-flags, they appear essentially as the usual flags.
The usual flags translate, at least in the transitive case, the Cartan
distribution on the jet space of a function of one variable. Multi-flags
translate, in the transitive case, the same situation in the jet space
of several functions of one variable.} \dots

Therefore, a kind of multi-dimensional prolongation of distributions was
badly needed. Va\-rious discussions around the results of \cite{Montgomery2001}
(existing then in a preprint form) and (drafts of)~\cite{Kumpera2002} remained
inconclusive until the formulation of a general prolongation scheme.

One obtains the def\/inition of {\it generalized Cartan prolongation} (gCp
for short, p.~159 in \cite{Mormul2004a}) by replacing in the def\/inition
from \cite{Bryant1993}: `rank-2' by `rank-$(m+1)$', `2-dimensional' by
`$(m+1)$-dimensional', and `$\P\R^1$' by `$\P\R^m$'. While $\P D$ stands,
as there, for the projectivization of the bundle $D \longrightarrow M$.

If $D$ is a rank-$(m+1)$ distribution on a manifold $M$, then, regarding~$D$ as
a vector bundle, its projectivization $\pi: \P D \longrightarrow M$ is a bundle
over $M$ whose typical f\/iber~$(\P D)_p$ is the space of 1-dimensional linear
subspaces of the $(m+1)$-dimensional vector space $D_p$. Thus, the f\/ibers of~$\P D$ are isomorphic to~$\P\R^m$ as projective $m$-spaces. There is a canonical
rank-$(m+1)$ distribution~$D^{(1)}$ on~$\P D$ def\/ined by setting $D^{(1)}_\xi =
(\pi')^{-1}(\xi)$ for each linear subspace $\xi \subset D_p$. This distribution~$D^{(1)}$ is called the (generalized) Cartan prolongation of~$D$.

Let us repeat that the prolonged distribution $D^{(1)}$ has the same
rank $m + 1$ as the initial distribution $D$, but it lives on a much bigger
manifold, having $m$ dimensions more than the initial manifold $M$. Similarly
as for the classical Cartan prolongation, immersed $D$-curves have canonical
lifts `upstairs' tangent to $D^{(1)}$. So it is clear what the local
generators of $D^{(1)}$ are.
For instance, one takes an immersed $D$-curve realizing any given horizontal
direction $\xi$ `downstairs', then takes the direction of its canonical lift,
and adds the $m$-dimensional kernel of the dif\/ferential, taken at that
point-direction $\xi$, of the projection $\pi$. Strictly speaking, a curve
realizing the direction $\xi$ is {\it not} necessary. It suf\/f\/ices to take
the horizontal vectors alone and lift them upstairs, although only relatively.
That is, {\it modulo} the kernel of $\pi'$. Having local generators
of $D$ -- like in Section~3.3 of~\cite{Mormul2004a} -- one is thus able
to `microlocally' write generators of~$D^{(1)}$. (At this moment one
already touches upon polynomial visualisations of the gCp's put forward
in \cite{Mormul2004a} and reiterated, for $m = 2$, in Section~\ref{EKRs}
of the present paper.)

 We note that certain ingredients (but ingredients only) of the
above def\/inition of the generalized Cartan prolongation were dispersed
in the literature, cf.\ Remark 1 in \cite{Mormul2004a} for more on that.

We intend now to recall a local structural theorem generalizing Cartan's
theorem from Section~\ref{IMT} which has geometrical applications, mainly
to special multi-f\/lags. Namely, the assumptions in Theorem~\ref{cbh} could
be rephrased by avoiding mentioning~$D_2$ and placing the foliation $\F$
in a new context. In fact, those assumptions easily implied that there
existed a (unique) line subdistribution $E$ of $D$ preserving $D_1$,
$[E,D_1] \subset D_1$. The foliation $\F$ was the integral of $E$.
Driven by the def\/inition of gCp's, we were going to replace a line
subdistribution of a rank-2 distribution in Theorem~\ref{cbh} by an
involutive rank-$m$ subdistribution of a rank-$(m+1)$ one (that is,
by its corank-1 involutive subdistribution).

\begin{theorem}[\cite{Mormul2004a}]\label{str}
Suppose D is a rank-$(m+1)$ distribution on a manifold $M^{s+m}$ such that
 {\rm a)}~$D_1$ is a rank-$(2m+1)$ distribution on~$M$, and {\rm b)}~there
exists a corank-$1$ involutive subdistribution $E \subset D$ that preserves
$D_1$,  $[E,D_1] \subset D_1$. Then D is locally equivalent to the
generalized Cartan prolongation $\bigl(D_1/E\bigr)^{(1)}$ of $D_1$ reduced
modulo~$E$ $($that lives on the quotient manifold $M/\F$ of  dimension~$s$,
where $\F$ is the local $m$-dimensional foliation in~$M$ defined by~$E)$.
\end{theorem}

Attention. $M/\F$ is to be understood only locally, to avoid topological
complications. Note that $\dim M = 2m + 1$, i.e., $s = m + 1$
is not excluded in this theorem.

It appears that distributions emerging as the outputs of several
applications of this theorem are precisely the jet bundles for maps
$\R \rightarrow \R^m$ {\it together} with the neighbouring distributions
pref\/igured by Kumpera. This should come as no surprise, for the gCp's
were tailored for Theorem~\ref{str}, which in turn was tailored for
the objects Kumpera and Rubin wrote about~-- especially in the f\/irst
version of~\cite{Kumpera2002} which was 60 + pages long. In short, it
is Theorem~\ref{str} which underlies the theory of special multi-f\/lags.
In particular it `makes possible' for Theorem~\ref{not} in
Section~\ref{threee} to hold true.

Aiming at completing now the discussion of the def\/inition of special
multi-f\/lags, we note

\begin{proposition}\label{uni}
Suppose that there is a distribution $D \subset TN$ of corank $m$ bigger
than~$1$, possessing an involutive corank-$1$ subdistribution~$E$, and such
that $[D,D] = TN$. Then, at each point $p \in N$, the value of~$E$
is described by all $($local$)$ $1$-forms $\alpha$ on $N$ such that
$\left.\bigl(\alpha\wedge d\omega\bigr)\right|_D = 0$ for all
$($local$)$ $1$-forms~$\omega$ annihilating~$D$.

 Moreover, the Cauchy-characteristic module $L(D)$ of~$D$ sits
then inside $E$ and is an involutive corank-$m$ subdistribution of~$E$.
\end{proposition}

This proposition is crucial for special multi-f\/lags and, hence, also for the
subject of the paper. It is proved in detail in Appendix~\ref{APP}.

\begin{remark}\label{remark1} (a) Whenever a family $\widetilde{E}$ of subspaces
of $D$, over points in $N$ and of dimensions \`{a} priori possibly depending
on those points, is being pointwisely described by the 1-forms $\alpha$ as
in Proposition \ref{uni} and $\widetilde{E}$ happens to be of constant
dimension, then, in \cite[p.~5]{Kumpera2002} it is called the {\it covariant
subdistribution} $\widehat{D}$ of $D$. So then $\widehat{D} = \widetilde{E}$.

(b) Technically, the authors of \cite{Kumpera2002} arrive at the
covariant object not directly, but via the so-called polar spaces of $D^\perp$,
included in~$T^*N/D^\perp$.  They continue only when the polar spaces are
of constant dimension, independently of a point. In the situation in
Proposition \ref{uni} that constant dimensionality turns out to be
automatic (see Appendix~\ref{APP}).
\end{remark}

\begin{corollary}\label{unique}
The involutive subdistribution $F \subset D^1$ from Definition~{\rm \ref{definition1}}
is unique and is nothing but the covariant subdistribution $\widehat{D^1}$.
It automatically contains $L(D^1)$ as its corank-$m$ subdistribution.
\end{corollary}

\begin{remark} \label{remark2} (a) Alternatively, one could assume
in $\star\star\star$  in Def\/inition~\ref{definition1} that the covariant
subdistribution of $D^1$ exists and is involutive. For, in view of Lemma~1
in \cite{Kumpera2002}, such a subdistribution is automatically of corank~1
in $D^1$; the hypotheses in that lemma are satisf\/ied as  ${\rm rk}\,
[D^1, D^1]/D^1  = m > 1$.

(b) Equivalently, using Tanaka's and Yamaguchi's terminology
\cite{Tanaka,Yamaguchi1982,Yamaguchi1983} (well anterior to \cite{Kumpera2002}),
one could stipulate in  $\star\star\star$  that the {\it symbol subdistribution}
of $D^1/L(D^1)$, which is automatically of corank~1 here, be involutive.
See also the detailed discussion of the symbol subdistribution on pages
28--30 in~\cite{Yamaguchi1983}.
\end{remark}

\subsection{Monsters for special multi-f\/lags}

As it has been explicitly stated in \cite{Mormul2004a} in Remark 3, every germ
of a distribution generating a special $m$-f\/lag of length $r$ can be found
within the $r$-fold generalized Cartan prolongation of the tangent bundle
to $\R^{m+1}$. This follows directly from Theorem \ref{str} coupled with
the original version of the def\/inition of special multi-f\/lags given in
section 3 of \cite{Mormul2004a} (equivalent to the present Def\/inition~\ref{definition1},
in which the Cauchy-characteristic subdistributions are not explicitly used).
Just like Goursat monster's coming into being was a direct consequence
of Theorem~\ref{cbh}. In the light of Theorem~\ref{str}, locally universal
objects in the theory of special multi-f\/lags are very natural.
In~\cite{Mormul2004a} they were abbreviated by MS$k$FM (from Monster Special
$k$-Flags Manifold), and they should now be written as MS$m$FM, $m$, not $k$,
standing now for the width.

In the recent paper \cite{Shibuya} the gCp is named `Rank 1 Prolongation'.
The result of $r$ consecutive gCp's applied to the tangent bundle to a manifold
$M$ of dimension $m + 1$, is called there an {\it $m$-flag of length~$r$} and is
denoted by $(P^r(M), C^r)$. (Strictly speaking, the distribution $C^r$
{\it generates} such a f\/lag.) In order not to multiply symbols, we will
adopt the notation $P^r(M)$ in the present paper, with a modest manifold
$M = \R^3$.

\section{Singularities of special 2-f\/lags}\label{threee}

It follows from the classical work \cite{Cartan1910} that special 2-f\/lags
of length~1 are homogeneous: they are identical around any point and
hence feature no singularities at all. Here are two examples of rank-3
distributions generating special 2-f\/lags of length 2. One of them is
still homogeneous and the other one has a singular locus of codimension 1.
The f\/irst example is the jet bundle on~$J^2(1,2)$,
\begin{equation}\label{ca_2}
\left(\frac{\p}{\p t} + x_1\frac{\p}{\p x} + y_1\frac{\p}{\p y} +
x_2\frac{\p}{\p x_1} + y_2\frac{\p}{\p y_1},  \frac{\p}{\p x_2},
  \frac{\p}{\p y_2}\right)
\end{equation}
(it is generalized to bigger lengths in Example~\ref{example1} below).
The second one is the following non-homogeneous object
\begin{equation}\label{ex_2}
\left(x_2\Bigl(\frac{\p}{\p t} + x_1\frac{\p}{\p x} +
y_1\frac{\p}{\p y}\Bigr) + \frac{\p}{\p x_1} + y_2\frac{\p}{\p y_1},
 \frac{\p}{\p x_2} ,  \frac{\p}{\p y_2}\right),
\end{equation}
which is singular on the hypersurface $\{x_2 = 0\}$. The fact that these
two distributions are non-equivalent as germs at $0 \in \R^7$ will be
clear in the next section. In fact, this is the starting point for the
theory proposed in the paper. In width and length both equal to two,
the object~(\ref{ca_2}) is the local model for the [generic and the
only one open] singularity class~1.1, while the object~(\ref{ex_2}) is
the model for the codimension-one singularity class~1.2 (see Sections~\ref{EKRs} and~\ref{refin} for precise def\/initions). These two classes
build up the stratif\/ication of germs of special 2-f\/lags when the length
$r = 2$. The objective of the paper is to do the same in any length.
Regarding further examples, a f\/iner instance of a couple of nonequivalent
2-f\/lags (of length 4) is given in Appendix~\ref{concr},
with the aim of illustrating the main constructions of the paper.

\subsection{Sandwich diagram for special 2-f\/lags}\label{SD}

Special multi-f\/lags, and in particular special 2-f\/lags, appear, from
one side, to be rich in singularities, and from the other, to possess
f\/inite-parameter families of local pseudo-normal forms, with no functional
moduli, constructed in~\cite{Mormul2004a}. The respective tree of normal
forms is very natural and emerges in a transparent way from the sequences
of gCp's being at work. Multi-parameter normal forms, in the case of
2-f\/lags dealt with in the present paper, are indexed by certain words
over the alphabet $\{\1,\2,\3\}$ of length equal to f\/lag's length.

On the other hand, a basic partitioning in the world of special
multi-f\/lags is a stratif\/ication into singularity classes proposed
in the preprint \cite{Mormul2003} and reproduced, for 2-f\/lags, below.
In their turn, the singularity classes for special 2-f\/lags are encoded
by certain words over the alphabet $\{1,2,3\}$ of length equal to
f\/lag's length.

Both partitions exist in their own rights, with no apparent
relation to each other. A f\/irst (modest) step towards throwing bridges
is the concept of sandwich classes (Section~\ref{san}), followed by
Corollary~\ref{relation} which makes use of that concept.

While the eventual aim of the paper, earlier undertaken in~\cite{Mormul2005} and interrupted, is to identify these two vocabularies:
to show that words over $\{\1,\2,\3\}$ and words over $\{1,2,3\}$
label precisely the same aggregates of germs of special 2-f\/lags -- see
Theorem \ref{2003}. A similar issue for multi-f\/lags of widths bigger
than~2 will be addressed in a future work.

Our initial requirements $\star\star$ and $\star\star\star$
are visualised best in a {\it sandwich diagram}\footnote{So-called
after a similar diagram assembled for Goursat distributions, or 1-f\/lags,
in \cite{Montgomery2001}.}
\[
 \begin{array}{@{}ccccccccccccc}
TM = D^0 & \supset & D^1 & \supset & D^2 & \supset & \cdots & \supset
& D^{r-1} & \supset & D^r & & \\
& & \cup & & \cup & & & & \cup & & \cup & & \\
& & F & \supset & L(D^1) & \supset & \cdots & \supset & L(D^{r-2})
& \supset & L(D^{r-1}) & \supset & L(D^r) = 0.
\end{array}
\]
The inclusions in its lower line are due to the Jacobi identity
($L(D^{j-1}) \supset L(D^j)$) and to Corollary~\ref{unique}
($F \supset L(D^1)$). All vertical inclusions in this diagram are of
codimension one while all drawn horizontal inclusions are of codimension~2.
The squares formed by these inclusions can be perceived as certain `sandwiches'.
For instance, in the utmost left sandwich $F$ and $D^2$ are as if f\/illings
while $D^1$ and $L(D^1)$ constitute the covers (of dimensions dif\/fering by~3).
At that, the sum (=\,3) of codimensions, in $D^1$, of $F$ and $D^2$ equals
the dimension of the quotient space $D^1/L(D^1)$, so that it is natural to
ask how the 2-dimensional {\it plane} $F/L(D^1)$ and the line $D^2/L(D^1)$
are mutually positioned in $D^1/L(D^1)$. Similar questions also arise
in further sandwiches `indexed' by the upper right `vertices'
$D^3,D^4,\dots,D^r$.

\subsection[Analogues for special 2-flags of Kumpera-Ruiz classes]{Analogues for special 2-f\/lags of Kumpera--Ruiz classes}
\label{san}

We thus f\/irst divide all existing germs of special 2-f\/lags of length
$r$ into $2^{r-1}$ pairwise disjoint {\it sandwich classes} depending
on the geometry of the distinguished spaces in the sandwiches (at the
reference point for a germ), and label those aggregates by words of length
$r$ over the alphabet~$\{$1,\underline{2}$\}$ starting (on the left)
with 1, having the second letter~\underline{2} if\/f $D^2(p) \subset F(p)$,
and for $3 \le j \le r$ having the $j$-th letter~\underline{2} if\/f
$D^j(p) \subset L(D^{j-2})(p)$.

It follows immediately from this def\/inition that the sandwich classes
are pairwise disjoint. On the other hand, it is not yet clear if each of
them is actually nonempty; this follows only from Corollary~\ref{relation}
below.

{\sloppy
The construction of sandwich classes points to possible non-transverse
situations in the sandwiches. For instance, the second letter in a sandwich
label is \underline{2} if\/f the line $D^2(p)/L(D^1)(p)$ is inclu\-ded in the plane
$F(p)/L(D^1)(p)$, both the line and plane sitting in the 3-space $D^1(p)/L(D^1)(p)$.
And it is similarly in further sandwiches. This resembles the Kumpera--Ruiz classes
of Goursat germs constructed in~\cite{Montgomery2001}. The number of sandwiches
in length $r$ then was $r-2$ (and so the $\#$ of~KR classes $2^{r-2}$) due to the
degenerate form of the covariant distribution of~$D^1$: $\widehat{D^1} = L(D^1)$.
Now, for 2-f\/lags this number is~$r-1$, because the covariant distribution of~$D^1$ dif\/fers from $L(D^1)$, and gives rise to one additional sandwich.

}

How can one establish if such virtually created sandwich classes really
materialize? And, if so, is it possible to sort them further?

 We shall produce a huge variety of polynomial germs at $0 \in \R^N$,
of rank-3 distributions, where $N$ will be odd and possibly be very large. It is
important that certain variables $x_j$ will appear in them in a shifted form
$b + x_j$, and it will always be an issue if such shifting constants are rigid
with respect to the local classif\/ication or subject to further simplif\/ications.
More precisely, for each $\k \in \{\1,\2,\3\}$ we are going to def\/ine
an {\it operation $\bf k$} producing new rank-3 distributions from previous ones.
Technically, its outcome (and especially the indices of new incoming variables)
will also depend on how many operations were done {\it before $\bf k$}.

The result of {\bf k}, being performed as an $l$-th operation in
a succession of operations, on a~distribution $(Z_1,Z_2,Z_3)$ def\/ined
in the vicinity of $0 \in \R^s(u_1,\dots,u_s)$, is a new rank-3 distribution~-- the germ at $0 \in \R^{s+2}(u_1,\dots,u_s,x_l,y_l)$, generated by
the vector f\/ields
\[
Z'_1 = \begin{cases}
Z_1 + (b_l + x_l)Z_2 + (c_l + y_l)Z_3 , & \text{when \ {\k} =
{\1}},\\
x_lZ_1 + Z_2 + (c_l + y_l)Z_3 , & \text{when \ {\bf k} = {\2}},\\
x_lZ_1 + y_lZ_2 + Z_3 , & \text{when \ {\k} = {\3}}\end{cases}
\]
and $Z'_2 = \frac{\p}{\p x_l}$, $Z'_3 = \frac{\p}{\p y_l}$.
 Here $b_l$ and/or $c_l$ are real parameters whose values
are specif\/ied later, when one applies these operations to concrete objects.
For any subsequent such operation (one will need to perform many of them)
it is important that these local generators are written precisely in this
order, yielding together a new `longer' or more involved distribution
$(Z'_1,Z'_2,Z'_3)$. Later (in Section~\ref{pf}) we will write more
compactly $X_l = b_l + x_l$, $Y_l = c_l + y_l$.

\subsection{Def\/inition of EKR's}\label{EKRs}

Extended Kumpera--Ruiz pseudo-normal forms (EKR for short), of length
$r \ge 1$, denoted by $\j_1.\j_2\dots\j_r$, where $\j_1,\dots,
\j_r \in \{\1,\2,\3\}$ and depending on numerous real parameters
within a~f\/ixed symbol $\j_1.\j_2\dots\j_r$, are def\/ined inductively,
starting from the distribution
\begin{equation}\label{zero}
\left(\frac{\p}{\p t},\frac{\p}{\p x_0},\frac{\p}{\p y_0}\right)
\end{equation}
understood in the vicinity of $0 \in \R^3(t,x_0,y_0)$;
this full tangent bundle to a 3-space is encoded by an empty label.
(The name `EKR' was coined in the work~\cite{Pasillas}, although the
very method of producing local visualisations of special multi-f\/lags
was not correct there. Namely, the authors of \cite{Pasillas} arrived
only at the operations~\1 and~\2. In fact, their relevant operations
are just~\1 and~\2 modulo reindexations in the $m$-tuples of their
variables $x^j_1,x^j_2,\dots,x^j_m$ ($j = 0,1,\dots,n$)
and similar reindexations in the $m$-tuples of their vector f\/ields
$\kappa^j_1,\kappa^j_2, \dots,\kappa^j_m$ ($j = 1,2,\dots,n$),
cf.~\cite[pp.~112--113]{Pasillas}.
While the operation $\3$ is necessary already for $m = 2$, as shows
Proposition~1(iv) in \cite{Mormul2004a} and the entire message of
the present article. Likewise, operation~{\bf 4} would turn out necessary
from width~3 and length~4 onwards, operation~{\bf 5} from width 4 and
length~5 on, etc.)  Assume that the family of pseudo-normal forms
$\j_1\dots\j_{r-1}$ is already constructed and written in coordinates
that go along with the operations: f\/irst $\j_1$, then $\j_2$ and so
on up to $\j_{r-1}$ (the distribution (\ref{zero}) when $r - 1 = 0$).
Then the normal forms subsumed under the symbol $\j_1\dots\j_{r-1}.
\j_r$ are the outcome of the operation $\j_r$ performed as the
operation number~$r$ over the distributions $\j_1\dots\j_{r-1}$.

For a moment, it is nearly directly visible that every EKR is
a special 2-f\/lag of length equal to the number of operations used
to produce it. In particular, it is easy to predict what the involutive
subdistributions of ranks $2,4,\dots,2r$ are; see also Proposition~\ref{posi}
below. The point is that locally the converse is also true.

\begin{theorem}\label{not}
Let a rank-$3$ distribution D generate a special $2$-flag of length
$r \ge 1$ on a mani\-fold~$M^{2r+3}$. For every point $p \in M$, the distribution~$D$ is equivalent in a neighbourhood of $p$ to a~certain EKR $\j_1.\j_2\dots\j_r$
in a neighbourhood of $0 \in \R^{2r+3}$, by a local diffeomorphism that
sends~$p$ to~$0$. Moreover, that EKR can be taken such that
$\j_1 = {\1}$ and the first letter~\2, if any, appears before
the first letter~\3 $($if any$)$.
\end{theorem}

This theorem is just the specif\/ication of Theorem~3 in \cite{Mormul2004a}
to the special 2-f\/lags. In particular, the restriction on EKR's codes
in it is the specif\/ication to the width $m = 2$ of the general rule of
the least upward jumps put forward in~\cite{Mormul2004a} and already
brief\/ly explained in Section~\ref{IMT}.

This rule looks modest in width 2. It becomes more and more
restrictive in larger widths $3,4,\dots$. Despite this, the idea
standing behind it is simple. At a new stage, one Cartan-prolongs in the
vicinity of a direction~$\xi$. What operation could one use for a local
description of that Cartan prolongation? Basically, any operation whose
pivot is not perpendicular to~$\xi$. Now suppose additionally that
{\it all} such operations have their numbers (or: indices) {\it higher}
than the indices of operations used before that stage. The rule under
discussion says that one should choose the operation which has the
{\it lowest} index among the not-yet-used indices. Technically,
it boils down to a {\it reindexation} of the `new' coordinates having
those higher indices. Then such a reindexation can safely be extended onto
the `old' coordinates bound to operations at earlier stages, not af\/fecting
the numbering of those earlier operations. Thus, inductively, one is able
to obey the rule of the least upward jumps. Details can be traced down
in~\cite[pp.~167--168]{Mormul2004a}.

We stress that possible constants in the EKRs representing a given
germ $D$ (i.e., the constants in the EKRs in Theorem~\ref{not})
are not, in general, def\/ined uniquely.

\begin{example}\label{example1} The EKR $\1.\1\dots\1$ ($r$ letters \1)
subsumes a vast family of dif\/ferent pseudo-normal forms -- germs at
$0 \in \R^{2r+3}$ parametrized by real parameters $b_1,c_1,\dots,b_r,c_r$.
Under a closer inspection (Theorem 1 in~\cite{Kumpera2002}), they all
are pairwise equivalent, and are equivalent to the classical jet bundle~-- Cartan distribution~-- on the space $J^r(1,2)$ of the $r$-jets
of functions $\R \rightarrow \R^2$, given by the Pfaf\/f\/ian equations
\[
dx_j - x_{j+1}dt = 0 = dy_j - y_{j+1}dt,\qquad
j = 0,1,\dots,r-1.
\]
All distribution germs in all other EKRs are not equivalent to
the jet bundles; this follows from Corollary~\ref{relation} below.

Let us note that the question of a geometric characterization of
Cartan distributions as such was addressed in many papers. In full
generality (for all jet spaces $J^r(n,m)$) that question was
answered only in 1983 in~\cite{Yamaguchi1983}.
\end{example}

\subsection{The EKR's versus sandwich classes}\label{versus}

What kind of a relationship does there exist between the sandwich
class of a~given germ of a special 2-f\/lag and its all possible EKR
presentations? In order to answer, we note

\begin{proposition}\label{posi}
If a distribution $D = D^r$ generating a special $2$-flag of length
$r \ge 1$ is presented in any EKR form on $\R^{2r+3}(t,x_0,
y_0,\dots,x_r,y_r)$, then the members of the associated subflag in
the sandwich diagram for $D^r$ are canonically positioned as follows
\begin{itemize}\itemsep=0pt
\item
$F = \bigl(\p/\p x_1,\p/\p y_1,\p/\p x_2,\p/\p y_2,\dots,
\p/\p x_r,\p/\p y_r\bigr)$,
\item
$L(D^j) = \bigl(\p/\p x_{j+1},\p/\p y_{j+1},\dots,\p/\p x_r,
\p/\p y_r\bigr)$ for $1 \le j \le r - 1$,
\item
$L(D^r) = (0)$.
\end{itemize}
\end{proposition}

Proof is almost immediate, because the inclusions $\supset$ (not yet
equalities) are clear in view of the construction of the EKR's, while
the dimensions of spaces on both sides of these inclusions always
coincide by the def\/inition of special 2-f\/lags.

(These extremely simplif\/ied descriptions of the members
of the associated subf\/lag are the analogues of similar descriptions
holding true for Goursat f\/lags viewed in Kumpera--Ruiz coordinates.)

Proposition~\ref{posi} has an important corollary. Namely,

\begin{corollary}\label{relation}
Each given sandwich class in length $r$ having label $\E$ is the
aggregate of all germs admitting EKR visualisations
of the forms $\j_1\dots\j_{r-1}.\j_r$ such that $\j_l = \1$
$\Longleftrightarrow$ the $l$-th letter in $\E$ is $1$, for
$l = 1,2,\dots,r$.
\end{corollary}

Therefore, the basic singular phenomena of the pointwise inclusions
in sandwiches do narrow (to~\2 and~\3) the pool of operations available
at the relevant steps of producing EKR visualisations for special 2-f\/lags.
The nonemptiness of sandwich classes follows.
Moreover, they are embedded submanifolds in the monster manifolds
$P^r(\R^3)$ of codimensions equal to the number of letters $\underline{2}$
in their codes. (We do not dwell on this any longer because by far more
important are smaller bricks, or singularity classes, building up
sandwich classes.)

\begin{proof} $\j_1$ is by default \1 and the f\/irst letter
in $\E$ is, by def\/inition, 1. Consider now $\j_l$, $l \ge 2$,
and recall that the operation $\j_l$ transforms certain EKR
$(Z_1,Z_2,Z_3)$ of length $l - 1$ into an EKR
$(Z'_1,Z'_2,Z'_3)$ of length $l$.  When $\j_l$ is
either \2 or \3 then, by def\/inition of these operations,
\begin{equation}\label{sanbis}
Z'_1 \equiv x_lZ_1 \ {\rm mod}\; (Z_2,Z_3),
\end{equation}
where $Z_2 = \frac{\p}{\p x_{l-1}}$ and $Z_3 = \frac{\p}{\p y_{l-1}}$.
(As for $Z'_2 = \frac{\p}{\p x_l}$ and $Z'_3 = \frac{\p}{\p y_l}$,
they cause no trouble in the discussion.) Whereas for $\j_l = \1$
we have $Z'_1 \equiv Z_1\ {\rm mod}\;(Z_2,Z_3)$ and the non-zero
vector $Z_1(0)$ is, by its recursive construction (in $l-1$ steps),
spanned by
\begin{equation}\label{KK}
\p/\p t,\ \p/\p x_0,\  \p/\p y_0,\ \dots,\ \p/\p x_{l-2},\
\p/\p y_{l-2}.
\end{equation}
Hence, in view of Proposition~\ref{posi}, $Z_1(0)$ does not lie in
$F(0)$ when $l = 2$, and in $L(D^{l-2})(0)$, when $l > 2$.
\end{proof}

\begin{remark}\label{remark3}
When $m = 1$ two operations, instead
of three (\1,\2,\3) in the present text, lead to the local
Kumpera--Ruiz pseudo-normal forms for Goursat f\/lags, evoked
already in Section~\ref{IMT}.
\end{remark}

\subsection{Singularity classes of special 2-f\/lags
ref\/ining the sandwich classes}\label{refin}

We recall from \cite{Mormul2003} how one passes from the sandwich classes
to {\it singularity classes}. In fact, to any germ $\F$ of a special
2-f\/lag we associate a word $\W(\F)$ over $\{$1,2,3$\}$, called a
`singularity class' of $\F$. It is a specif\/ication of the word `sandwich
class' for $\F$ (a word, recalling, over $\{$1,\underline{2}$\}$)
with the letters~\underline{2} replaced either by 2 or 3, depending
on the geometry of $\F$. It will be momentarily clear from the def\/inition
that $\W(\cdot)$ is an invariant of the local classif\/ication of f\/lags
with respect to dif\/feomorphisms in the base manifold.

Alternatively, if one restricts oneself to the locally universal
f\/lags of distributions $C^r$ living on~$P^r(\R^3)$, then $\W$ becomes
essentially a function of a point in $P^r(\R^3)$, and it will turn out
to be an invariant of the (local) symmetries of~$C^r$. That is, an invariant
of the local dif\/feomorphisms of~$\R^3$, inducing after $r$ prolongations
the symmetries of~$C^r$ on~$P^r(\R^3)$.

In the def\/inition that follows we keep the germ of a rank-3 distribution $D$
generating a special 2-f\/lag $\F$ of length $r$ on $M$ f\/ixed.

Suppose that in the sandwich class $\E$ of $D$ at $p$ there
appears somewhere, for the f\/irst time when going from the left, the letter
$\underline{2} = j_f$ ($j_f$ is assuredly not the f\/irst letter in $\E$)
and that there are in $\E$ other letters $\underline{2} = j_s$,
$f < s$, as well. We will specify each such $j_s$ to either 2 or 3.
(The specif\/ication of the f\/irst $j_f$ will be made later and will be easy.)
Let the nearest \underline{2} standing to the left to $j_s$ be $\underline{2}
= j_\nu$, $f \le \nu < s$. These two `neighbouring' letters~\underline{2}
are separated in $\E$ by $l = s - \nu - 1 \ge 0$ letters 1.

 The core of the construction consists in taking the
{\it small flag} of f\/lag's member $D^s$,
\[
D^s = V_1 \subset V_2 \subset V_3 \subset V_4 \subset V_5 \subset
\cdots,
\]
$V_{i+1} = V_i + [D^s,V_i]$, and then focusing on this new f\/lag's member
$V_{2l+3}$. Recall that, in the $\nu$-th sandwich, there holds the inclusion:
$F(p) \supset D^2(p)$, when $\nu = 2$, or else $L(D^{\nu-2})(p) \supset
D^\nu(p)$, when $\nu > 2$. This is a preparation to an important, turning
point decision.

 Namely, writing $V_{2l+3}(p)$ instead of $D^\nu(p)$ in the relevant
inclusion, and always controlling whether $\nu = 2$ or $\nu > 2$, means
specifying $j_s$ to 3. That is to say, $j_s = \underline{2}$ is being
specif\/ied to~3 if and only if $F(p) \supset V_{2l+3}(p)$ (when
$\nu = 2$) or else when $L(D^{\nu-2})(p) \supset V_{2l+3}(p)$ (when
$\nu > 2$) holds.

In this way all non-f\/irst letters \underline{2} in $\C$ are, one independently
of another, specif\/ied to~2 or~3. Having done that, one simply replaces the f\/irst
letter \underline{2} by 2, and altogether obtains a word over $\{1,2,3\}$.
It is the singularity class $\W(\F)$ of $\F$ at $p$. The word created by
such a procedure clearly satisf\/ies the least upward jumps rule.

This is the singularity class of a given 2-f\/lag at a point. So what is
an abstract singularity class in length $r$, what subset of the monster manifold
$P^r(\R^3)$ does it form? It is the union of all points in $P^r\bigl(\R^3\bigr)$
at which the {\it universal} f\/lag has a f\/ixed singularity class -- a f\/ixed word
of length $r$ over $\{1,2,3\}$ obeying the rule of least upward jumps. Hence
there emerges a partition of $P^r\bigl(\R^3\bigr)$ into abstract, pairwise disjoint
singularity classes.

\begin{example}\label{example2} In length 4 there exist (or: $P^4(\R^3)$ is
partitioned into) the following fourteen singularity classes: 1.1.1.1,
1.1.1.2, 1.1.2.1, 1.1.2.2, 1.1.2.3, 1.2.1.1, 1.2.1.2,
1.2.1.3, 1.2.2.1, 1.2.2.2, 1.2.2.3, 1.2.3.1, 1.2.3.2,
1.2.3.3.\footnote{In widths $\ge3$ the class 1.2.3.4
will show up as well, cf.\ Remark~\ref{remark4}(b).}

(Reiterating, the emptiness of certain singularity classes
has not been {\`a} priori excluded. Only Theorem \ref{2003} shows
that all singularity classes are not empty -- see the paragraph
after Corollary~\ref{lass}.)
\end{example}

 How many singularity classes do there exist for special
2-f\/lags of f\/ixed length?

\begin{proposition}\label{how}
The number of different singularity classes of special $2$-flags
of length $r \ge 3$ is $2 + 3 + 3^2 + \cdots + 3^{r-2} =
\frac{1}{2}(1 + 3^{r-1})$.
\end{proposition}

\begin{proof}
Let us recall that the class' code $j_1.j_2\dots j_r$ is
subject to the least upward jumps rule. Either it is 1.1$\dots$1,
or else it has the f\/irst from left letter $j_f = 2$ at the $f$-th
position, $2 \le f \le r$. For $f = r$ one gets just 1 class.
For $f = r-1$ the number of classes' codes is $3^1$, for $f = r-2$
that number is $3^2$, and so on downwards to $f = 2$, with the
respective number of such classes~$3^{r-2}$.
\end{proof}

\begin{remark}\label{remark4} (a) The singularity classes discussed in
the present paper are just the visible part of an iceberg. Their
counterparts in the Goursat world, the KR classes, are nothing but vague
approximations to the orbits of the local classif\/ication. Much f\/iner
are then {\it geometric classes} emanating from Jean's benchmark
contribution \cite{Jean} and otherwise pref\/igured in~\cite{Montgomery2001}.
They are described in detail in~\cite{Mormul2004b}\footnote{In a
dif\/ferent language using extensively classical Cartan prolongations
and projections in the Goursat monster tower, the geometric classes,
under the name of `RVT classes', have been recently very originally
treated in~\cite{Montgomery2008}.}. Although they, too, are encoded by
certain words over a three letters' alphabet, one should {\it by no means}
confuse them with singularity classes for special 2-f\/lags.
The question of further partitioning of singularity classes
for special 2-f\/lags (and/or generally for special multi-f\/lags)
is under investigation, if still open for the most part.

(b) Reiterating after Section~\ref{SD}, singularity classes
for {\it all} widths $m$ have been def\/ined in~\cite{Mormul2003}.
To give an idea of their numbers, let us, for example, f\/ix the length
$r = 7$. Then the numbers of dif\/ferent singularity classes of special
$m$-f\/lags, for $m \in \{1,2,\dots,6\}$, are as follows (for $m = 1$
counted are the KR classes):
\begin{center}
\begin{tabular}{|c|||c||c|c|c|c|c|}\hline
$m$ & $1$ & $2$ & $3$ & $4$ & $5$ & $6$\\ \hline
                \hline
$\#$ & $32$ & $365$ & $715$ & $855$ & $876$ & $877$ \\ \hline
\end{tabular}
\end{center}
The value 365 is the value for $r = 7$ of the expression given
in Proposition~\ref{how}.
\end{remark}

\begin{remark}\label{remark5} Theorem \ref{2003} naturally generalizes
to {\it wider} special f\/lags. For $m > 2$ the f\/irst ref\/inement
of a sandwich class -- a word over $\{$1,2,\underline{3}$\}$
(see~\cite{Mormul2003} for details) -- is not yet a singularity class.
But it is a purely geometric notion, imposing natural restrictions on
the EKRs representing germs that have a f\/ixed word $j_1.j_2\dots j_r$
over $\{$1,2,\underline{3}$\}$ (satisfying the least upward jumps
rule). If $\k_1.\k_2\dots \k_r$ is any such an EKR, then
$j_l = \min (k_l,\underline{3})$ for $l = 1,2,\dots,r$.
That is, $j_l = k_l$ for those $\k_l$'s that are equal to~\1 or~\2,
and $j_l = \underline{3}$ for all the remaining $\k_l$'s.

A proof of this generalization of Theorem \ref{2003} is
only technically more complex, but not more dif\/f\/icult than the one
presented in the following chapter.
\end{remark}

Last but not least, there arises a question concerning the
materializations of singularity classes for {\it concrete}
special 2-f\/lags. In fact, on each manifold $M$ of dimension
$2r + 3$, $r \ge 1$, bearing a~special 2-f\/lag of length $r$,
the shadows of universal singularity classes in $P^r(\R^3)$
always form~-- and not only for `generic' f\/lags~-- a very neat
stratif\/ication by embedded submanifolds whose codimensions
are directly computable. Namely, we have the following

\begin{proposition}\label{twice}
The codimension of an embedded in $M$ submanifold of the realization
of any fixed singularity class $\C$, if only nonempty, is equal to
\begin{equation}
{\rm the\ number\ of\ letters\ 2\ in\ }\C \,\,+
\,\,{\rm twice\ the\ number\ of\ letters\ 3\ in\ }\C.\tag{$*$}
\end{equation}
In particular, the same formula $(*)$ holds for each singularity
class $\C \subset P^r(\R^3)$. In this case $\C$ is automatically
nonempty because of the universality property of $P^r(\R^3)$: $\C$
is mapped by the relevant EKR coordinates into certain
$\R^{(r+1)2+1}$ bearing the EKR forms with the label identical
to the label of $\C$. Speaking differently, the monster manifold
$P^r(\R^3)$ carries a universal $($in length~$r)$ stratification
into nonempty singularity classes.
\end{proposition}

A sketched proof of this is postponed until after the proof of
Theorem \ref{2003} (Appendix~\ref{apen}). In turn,
once the codimensions are made explicit, another natural question
is that about the adjacencies existing among these classes. We do
not have a full answer to this question yet. We only know that,
in any f\/ixed length $r \ge 1$,

\begin{proposition}
The generic class  {\rm 1.1$\dots$1} is not adjacent to any other
singularity class. An adjacency $j_1.j_2\dots j_l\dots j_r \rightarrow
j_1.j_2\dots(j_l-1)\dots j_r$, $2 \le l \le r$, holds whenever
$j_l = 3$ or $j_l = 2$, provided, in the latter case, there is no
letter~$3$ past~$j_l$ $($i.e., among $j_{l+1},\dots,j_r)$.
\end{proposition}

For instance, $1.2.3 \rightarrow 1.2.2 \rightarrow 1.1.2 \rightarrow
1.1.1$, or else $1.2.3.2 \rightarrow 1.2.3.1 \rightarrow 1.2.2.1
\rightarrow \cdots$. To completely answer the question, a deep analysis
of EKRs (i.e., the ef\/fective realizations, or visualisations, of
special f\/lags) is needed.

\section{Proof of Theorem~\ref{2003}}\label{pf}

We assume that the reader remembers the way the sandwich classes
were ref\/ined to singularity classes in Section~\ref{refin}. In the
proof of Theorem \ref{2003} we stay within that same framework (and
notation) and assume that:
\begin{enumerate}\itemsep=0pt
\item[--] the $\nu$-th letter $j_\nu$ in $\C$ is {\it not}~1,

\item[--] there follow $l \ge 0$ letters 1 past $j_\nu$,

\item[--] the following letter $j_s$ is {\it not} 1, where $s = \nu+l+1$.
\end{enumerate}

Having $D = D^r$ in a not-yet-specif\/ied EKR form $\k_1.\k_2\dots \k_r$
we know by Corollary~\ref{relation} that $\k_\nu \ne \1$, $\k_{\nu+1}
= \cdots = \k_{s-1} = \1$, $\k_s \ne \1$. And we aim to show that
\begin{equation}
\k_s = \3\textrm{\quad if and only if} \ \ j_s = 3.\tag{$**$}
\end{equation}
Only this is an issue. For, the f\/irst from the left letter $\k_f \ne \1$
(if any) is \2 by the least upward jumps rule satisf\/ied by the labels of
EKR classes, and the corresponding letter $j_f$ in $\C$ is~-- by the
same Corollary~\ref{relation}~-- the f\/irst from the left not 1 letter
in~$\C$. Hence it is~2 by the very def\/inition of singularity classes.

As for $(**)$, in Section~\ref{light} we will show
that $\k_s = \2$ implies $j_s = 2$, and in section \ref{hard}
that $\k_s = \3$ implies $j_s = 3$. That will do, because $\k_s
\in \{\2,\3\} \Longleftrightarrow j_s \in \{2,3\}$ by
Corollary~\ref{relation}.

Prior to concrete computations, note that, automatically, the rank-3
distribution $D^s/L(D^s)$, generating a special 2-f\/lag of length $s$,
is in an EKR form $\k_1.\k_2\dots \k_s$. In the (rather long)
computations that follow we {\it skip} writing down this factoring
out by the Cauchy characteristics~$L(D^s)$. That is, we simply leave out
the variables with indices from $s + 1$ onwards, upon which~$D^s$ does
not depend (Proposition~\ref{posi}). Also, for space reasons, from now
on we shall just write~$\p_x$ and $\p_{x_k}$ instead of~$\p/\p x$ and
$\p/\p x_k$, respectively.

\subsection[Easier part: $\k_s = \2$]{Easier part: $\boldsymbol{\k_s = \2}$}
\label{light}

We f\/irst deal with the case $\k_s = \2$ and aim at showing that then
$j_s = 2$ (meaning the {\it non-}inclusion of $V_{2l+3}(0)$ in the
relevant member of the Cauchy-characteristic subf\/lag).

\begin{proof} Let us expand the f\/irst member of the small f\/lag of $D^s$
\begin{gather}
D^s = V_1=\Biggl(x_s\Bigl(\underline{x_\nu Z + \ast\bigl(\p_{x_{\nu-1}},\p_{y_{\nu-1}}
\bigr)} + \sum_{k = \nu}^{s - 2} \bigl(X_{k+1}\p_{x_k} + Y_{k+1}\p_{y_k}\bigr)\Bigr)\nonumber\\
\phantom{D^s = V_1=\Biggl(}{}
+ \p_{x_{s-1}} + Y_s\p_{y_{s-1}}, \p_{x_s},\p_{y_s}\Biggr),\label{nu2}
\end{gather}
where the underlined summand is the leading generator of f\/lag's member $D^\nu$.
That is,
\[
D^\nu/L(D^\nu) = \Bigl(x_\nu Z + \ast\bigl(\p_{x_{\nu-1}},\p_
{y_{\nu-1}}\bigr),\p_{x_\nu},\p_{y_\nu}\Bigr)
\]
and the functions
$\ast$ depend on whether $\k_\nu$ is~\2 or~\3.
The capital letters $X$ and $Y$ stand, as in the end of Section~\ref{san},
for variables shifted by constants: $X_{k+1}=b_{k+1}+x_{k+1}$,  $Y_{k+1}
= c_{k+1} + y_{k+1}$.  By means of a straightforward step by step
computation, stopping at each {\it odd} member of the small f\/lag,
one shows that
\begin{gather}
\p_{x_{s-2}} + Y_s\p_{y_{s-2}} \in V_3,\nonumber\\
\p_{x_{s-3}} + Y_s\p_{y_{s-3}} \in V_5,
\nonumber\\
\cdots\cdots\cdots\cdots\cdots\cdots\nonumber\\
\p_{x_\nu} + Y_s\p_{y_\nu} \in V_{2l+1},
\nonumber\\
\label{shy}
Z + Y_s\bigl(\p_{x_{\nu - 1}},\p_{y_{\nu - 1}}\bigr) \in V_{2l+3},
\end{gather}
where $(\p_{x_{\nu - 1}},\p_{y_{\nu - 1}})$ stands for certain
combination of the versors $\p_{x_{\nu - 1}}$ and $\p_{y_{\nu - 1}}$.
In view of Proposition~\ref{posi}, these versors lie in $L(D^{\nu-2})$,
or in $F$ when $\nu = 2$. Whereas $Z(0)$ is, exactly as in Section~\ref{versus}, a nonzero combination of versors (\ref{KK}) for $l = \nu$
and as such sticks out of $L(D^{\nu - 2})(0)$, or of $F(0)$ when $\nu=2$.
Therefore~(\ref{shy}) alone implies that $V_{2l+3}(0)$ is not included
in $L(D^{\nu - 2})(0)$, or in $F(0)$ when $\nu = 2$. That is, $j_s = 2$.
\end{proof}

\subsection[Harder part: $\k_s = \3$]{Harder part: $\boldsymbol{\k_s = \3}$}\label{hard}

One should justify that now $j_s = 3$. That is, that there holds
the inclusion
\begin{equation}\label{hol}
V_{2l+3}(0) \subset \left(\p_{x_{\nu - 1}},  \p_{y_{\nu - 1}}, \ \p_
{x_\nu}, \p_{y_\nu}, \dots,  \p_{x_s},  \p_{y_s}\right).
\end{equation}

\begin{proof} The initial object $D^s = V_1$ is now dif\/ferent from
(\ref{nu2}).  Namely,
\begin{gather*}
V_1  =\Biggl(\!x_s\Bigl(x_\nu Z + \ast\bigl(\p_{x_{\nu-1}}, \p_{y_{\nu-1}}
\bigr) + \sum_{k = \nu}^{s - 2} \bigl(X_{k+1}\p_{x_k}
 + Y_{k+1}\p_{y_k}
\bigr)\Bigr) + \underline{y_s\p_{x_{s-1}} +  \p_{y_{s-1}}} ,
\p_{x_s} ,  \p_{y_s}\!\Biggr),\!
\end{gather*}
where $\ast$ stands for certain functions depending on the value
of~$\k_\nu$. Note the only dif\/ference, in the underlined part, with
the leading generator in~(\ref{nu2}). This slight dif\/ference will turn
out to be decisive in the output $V_{2l+3}$. Let us compute carefully
some f\/irst members of the small f\/lag of~$D^s$:
\begin{gather}
V_2  =  \Biggl(x_s\Bigl(x_\nu Z + \ast\bigl(\p_{x_{\nu-1}},
\p_{y_{\nu-1}}\bigr) +  \sum_{k = \nu}^{s - 2} \bigl(X_{k+1}\p_{x_k}
+ Y_{k+1}\p_{y_k}\bigr)\Bigr),  \p_{x_{s-1}} ,
 \p_{y_{s-1}}, \p_{x_s},  \p_{y_s}\Biggr),\nonumber\\
V_3  =  \left(V_2, x_s\p_{x_{s-2}}, x_s\p_{y_{s-2}}, y_s\p_{x_{s-2}} + \p_{y_{s-2}}\right),
\label{y-s}\\
V_4  =  \Bigl(V_3 ,  \p_{x_{s-2}}, \p_{y_{s-2}}, x_s^{2}
\p_{x_{s-3}}, x_s^{2}\p_{y_{s-3}}, x_s\bigl(y_s\p_{x_{s-2}}
+  \p_{y_{s-2}}\bigr)\Bigr).\nonumber
\end{gather}
Acting likewise, one keeps expressing $V_{n+1}$ by the previous module
$V_n$ and a set of simple vector f\/ield's generators, of the cardinality
growing linearly with~$n$, for $n \le l+2$. The modules $V_{l+2}$ and
$V_{l+3}$ are the most important in this process of computing.
The reader will see that in $V_{l+3}$ for the f\/irst time there appears
the f\/ield $Z$ standing alone, only with a monomial factor of high degree.
That f\/ield requires a particular care; in the situation $\k_s = \2$
it has been responsible for the failure of the inclusion. Strictly
speaking, the modules $V_{l+2}$ and $V_{l+3}$ look dif\/ferently
depending on the parity of~$l$:
\begin{enumerate}\itemsep=0pt
\item[$\star$] $l = 2k - 1$,  $k \ge 1$,  or else

\item[$\star\star$] $l = 2k$,  $k \ge 0$.
\end{enumerate}

 However, these dif\/ferences are not fundamental and one common
technique works in both situations. But a choice {\it is} necessary
when it comes to details. So for the presentation in the text we
choose $\star$.

For $l$ odd, $V_{l+2}$ is the module generated by $V_{l+1}$
and by the following set of generators:
\begin{itemize}\itemsep=0pt
\item
$x_s\p_{x_{\nu+k-1}} ,\  x_s\p_{y_{\nu+k-1}} ,\
y_s\p_{x_{\nu+k-1}} +  \p_{y_{\nu+k-1}}$;
\item
$x_s^{3}\p_{x_{\nu+k-2}} ,\  x_s^{3}\p_{y_{\nu+k-2}},
\  x_s^{2}\left(y_s\p_{x_{\nu+k-2}} +\p_{y_{\nu+k-2}}
\right)$;
\item
 $\cdots\cdots\cdots\cdots\cdots\cdots\cdots\cdots\cdots\cdots\cdots\cdots\cdots\cdots\cdots$
\item
$x_s^{2k-1}\p_{x_\nu},\ x_s^{2k-1}\p_{y_\nu} ,
\ x_s^{2k-2}\left(y_s\p_{x_\nu} +\p_{y_\nu}\right)$.
\end{itemize}
It has been straightforward to see that $V_{l+1}(0)$ is included in
the RHS of~(\ref{hol}). Hence so is $V_{l+2}(0)$. In turn, $V_{l+3}$
is the previous module $V_{l+2}$ extended by the generators
\begin{itemize}\itemsep=0pt
\item
$\p_{x_{\nu+k-1}},\ \p_{y_{\nu+k-1}}$;
\item
$x_s^{2}\p_{x_{\nu+k-2}},\  x_s^{2}\p_{y_{\nu+k-2}},
\  x_s\left(y_s\p_{x_{\nu+k-2}} + \p_{y_{\nu+k-2}}\right)$;
\item
 $\cdots\cdots\cdots\cdots\cdots\cdots\cdots\cdots\cdots\cdots\cdots\cdots\cdots\cdots\cdots\cdots$
\item
$x_s^{2k-2}\p_{x_\nu},\  x_s^{2k-2}\p_{y_\nu},
\  x_s^{2k-3}\left(y_s\p_{x_\nu} +  \p_{y_\nu}\right)$;
\item
$x_s^{l+1}Z,\  x_s^{l+1}\bigl(\p_{x_{\nu-1}},
\p_{y_{\nu-1}}\bigr),\ x_s^{l}y_sZ +
x_s^{l}\bigl(\p_{x_{\nu-1}},\p_{y_{\nu-1}}\bigr)$
\end{itemize}
(remember that $2k = l+1$). This gives that $V_{l+3}(0)$ is included
in the RHS of~(\ref{hol}). Having $V_{l+3}$ thus described is a turning
point in the proof. Indeed, there remains exactly~$l$ steps from~$V_{l+3}$
to $V_{2l+3}$. Because of that, in the bottom line of the new generators
in $V_{l+3}$, all terms with degree $l+1$ monomials are irrelevant as they
give rise only to terms vanishing at~0 under~$l$ Lie multiplications.

The remaining terms $x_s^{l}\bigl(\p_{x_{\nu - 1}},
\p_{y_{\nu - 1}}\bigr)$ in that bottom line could contribute at 0
only by means of dif\/ferentiating that monomial $l$ times during
$l$ Lie bracketings still to be performed, because of a degree $l$
monomial in them. Consequently they yield only the output sitting (at~0)
in $\bigl(\p_{x_{\nu - 1}}, \p_{y_{\nu - 1}}\bigr)(0)$.

All in all, in $l$ steps, the bottommost line of generators of $V_{l+3}$ will give
rise uniquely to vector f\/ields having at 0 values sitting in the RHS of~(\ref{hol}).
Therefore only the remaining lines of generators are relevant for the answer
to the question whether~(\ref{hol}) holds.
Thus, for answering this question, $V_{l+3}$ can be {\it replaced}
till the end of computations by the module $\overline{V}_{l+3}$
generated by $V_{l+2}$ and the smaller set of vector f\/ields
\begin{itemize}\itemsep=0pt
\item
$\p_{x_{\nu+k-1}} ,\  \p_{y_{\nu+k-1}}$;
\item
$x_s^{2}\p_{x_{\nu+k-2}} ,\  x_s^{2}\p_{y_{\nu+k-2}} ,
\  x_s\left(y_s\p_{x_{\nu+k-2}} +  \p_{y_{\nu+k-2}}\right)$ ;
\item
$\cdots\cdots\cdots\cdots\cdots\cdots\cdots\cdots\cdots\cdots\cdots\cdots\cdots\cdots\cdots\cdots$
\item
$x_s^{2k - 2}\p_{x_\nu},\ x_s^{2k - 2}\p_{y_\nu},
\  x_s^{2k-3}\left(y_s\p_{x_\nu} + \p_{y_\nu}\right)$.
\end{itemize}
So $\overline{V}_{l+3}$ is to be Lie bracketed $l$ times with~$V_1$.
Then the result taken at~0, $\overline{V}_{\!2l+3}(0)$, is to be checked
for its inclusion in the RHS of~(\ref{hol}), and that {\it would} f\/inish
the proof. Yet, as it stands, it is not transparent at all, and a series
of further simplif\/ications is needed. Computing now the next module
$\overline{V}_{l+4} = \overline{V}_{l+3} + \big[V_1,\overline{V}_{l+3}\big]$, one sees that it is the module generated
by $\overline{V}_{l+3}$ and the following collection of vector f\/ields
\begin{itemize}\itemsep=0pt
\item
$x_s\p_{x_{\nu+k-2}},\  x_s\p_{y_{\nu+k-2}} ,\
y_s\p_{x_{\nu+k-2}} +  \p_{y_{\nu+k-2}}$;
\item
$x_s^{3}\p_{x_{\nu+k-3}},\ x_s^{3}\p_{y_{\nu+k-3}},
\ x_s^{2}\left(y_s\p_{x_{\nu+k-3}} +\p_{y_{\nu+k-3}}
\right)$;
\item
$\cdots\cdots\cdots\cdots\cdots\cdots\cdots\cdots\cdots\cdots\cdots\cdots\cdots\cdots\cdots\cdots$
\item
$x_s^{2k-3}\p_{x_\nu},\ x_s^{2k-3}\p_{y_\nu},
\ x_s^{2k-4}\left(y_s\p_{x_\nu} + \p_{y_\nu}\right)$;
\item
$x_s^{l}Z,\ x_s^{l}\bigl(\p_{x_{\nu-1}},
\p_{y_{\nu-1}}\bigr),\ x_s^{l-1}y_sZ +
x_s^{l-1}\bigl(\p_{x_{\nu-1}}, \p_{y_{\nu-1}}\bigr)$
\end{itemize}
(remember that $2k-1 = l$). Hence $\overline{V}_{l+4}(0)$
sits in the RHS of (\ref{hol}).

 From now on the arguments start to repeat themselves.
Only $l-1$ Lie bracketings with $V_1$ remain to be done, hence the
last line of new generators for $\overline{V}_{l+4}$ is irrelevant
for the sought inclusion of $\overline{V}_{2l+3}(0)$ in the RHS
of~(\ref{hol}). Consequently, $\overline{V}_{l+4}$ can be replaced
until the end of computations by the module $\overline{\overline{V}}_
{l+4}$ generated by $\overline{V}_{l+3}$ and the smaller set of
vector f\/ields
\begin{itemize}\itemsep=0pt
\item
$x_s\p_{x_{\nu+k-2}},\ x_s\p_{y_{\nu+k-2}},\
y_s\p_{x_{\nu+k-2}} + \p_{y_{\nu+k-2}}$;
\item
$x_s^{3}\p_{x_{\nu+k-3}},\ x_s^{3}\p_{y_{\nu+k-3}},
\quad x_s^{2}\left(y_s\p_{x_{\nu+k-3}} + \p_{y_{\nu+k-3}}
\right)$;
\item
$\cdots\cdots\cdots\cdots\cdots\cdots\cdots\cdots\cdots\cdots\cdots\cdots\cdots\cdots\cdots\cdots$
\item
$x_s^{2k-3}\p_{x_\nu},\ x_s^{2k-3}\p_{y_\nu},
\quad x_s^{2k-4}\left(y_s\p_{x_\nu} + \p_{y_\nu}\right)$.
\end{itemize}
In its turn, the module $\overline{\overline{V}}_{l+4}$ will lead in
$l - 1$ steps to a module $\overline{\overline{V}}_{2l+3}$ which will
suf\/f\/ice for the verif\/ication of the inclusion as well. In the f\/irst of these
steps one is to:
\begin{enumerate}\itemsep=0pt
\item[--] compute the module $\overline{\overline{V}}_{l+5}
= \overline{\overline{V}}_{l+4} + \big[V_1,\overline{\overline{V}}_
{l+4}\big]$,

\item[--] check its inclusion at 0 in the RHS of (\ref{hol}), and then

\item[--] leave out its bottommost line of new generators, irrelevant
for the verif\/ication of the inclusion of $\overline{\overline{V}}_{2l+3}(0)$
in the RHS of (\ref{hol}).
\end{enumerate}
And then to proceed similarly in the remaining steps.

Summarizing, the critical part of the procedure applied in the proof
of Theorem~\ref{2003} consists of the following bipartite steps, having
numbers $l + 3 + \tau$, where $\tau \in \{1,\dots,l-2\}$.

 {\it Firstly} in checking the inclusion at 0, of the due reduction
of $V_{l + 3 + \tau}$, in the RHS of~(\ref{hol}).
{\it Secondly} in deleting the most involved, and not important
for the eventual inclusion of $V_{2l + 3}(0)$, among the {\it new} vector
f\/ield's generators showing up in the Lie product of $V_1$ with the mentioned
reduction of $V_{l + 3 + \tau}$.

 After reducing the problem in $l - 2$ steps in the described way,
the situation is as follows (we use the $\ \widetilde{ }\ $ symbols instead
of writing many bars):

One knows already from the step number $l + 3 + l - 2$ that
a)~the reduced module $\widetilde{V_{2l+1}}$ is such that $\widetilde
{V_{2l+1}}(0)$ is included in the RHS of (\ref{hol}), and b)~the newly
emerging module $\widetilde{\widetilde{V_{2l+2}}}$, again suf\/f\/icient
for the verif\/ication of the inclusion in question because only irrelevant
(for that verif\/ication) new generators have {\it just} been deleted,
is generated by $\widetilde{V_{2l+1}}$ and by a tiny set of generators
\begin{itemize}\itemsep=0pt
\item
$\p_{x_\nu} ,\  \p_{y_\nu}$.
\end{itemize}
Hence, at 0, it is included in the RHS of (\ref{hol}) as well. Moreover,
the module
\[
\widetilde{\widetilde{V_{2l+3}}}  =  \widetilde{\widetilde{V_{2l+2}}}
+ \big[V_1,\widetilde{\widetilde{V_{2l+2}}}\big],
\]
when evaluated at 0, is included in the RHS of~(\ref{hol}), too.
Indeed, $\p_{y_\nu}$ bracketed with $V_1$ gives the f\/ields in
$\bigl(\p_{x_{\nu-1}},\p_{y_{\nu-1}}\bigr)$, and $\p_{x_\nu}$
bracketed with the generators of $V_1$ gives just $x_sZ$ which
vanishes at~0.

But, by virtue of the adopted procedure of reductions, the inclusion
of $\widetilde{\widetilde{V_{2l+3}}}(0)$ in the RHS of~(\ref{hol})
is equivalent to the similar inclusion of $V_{2l+3}(0)$. Therefore
(\ref{hol}) holds. That is, by Proposition \ref{posi}, $V_{2l+3}(0)
\subset L(D^{\nu-2})(0)$ when $\nu > 2$, and $V_{2l+3}(0) \subset F(0)$
when $\nu = 2$. That is to say, $j_s = 3$.

The reasoning in the even situation $\star\star$ is analogous,
with indices and exponents interweaving with those of the odd case $\star$.
Theorem \ref{2003} is now fully proved.
\end{proof}

\appendix

\section{Proof of Proposition \ref{twice}}\label{apen}

For a given distribution $D$ on a manifold $M$ we take any singularity
class $\C = j_1.j_2\dots j_r$ hit by (or: nonempty for)~$D$. There
can be many classes hit by $D$, if not necessarily all $\frac{1}{2}
\bigl(1 + 3^{r-1}\bigr)$ existing in the length $r$ under discussion.
Then take any point~$p$ of the relevant singularity locus~$S$.
Here is an argument that~$S$ around $p$ is an embedded submanifold
of codimension stated in the proposition.

\begin{proof}
We take any f\/ixed polynomial presentation for $D$ around $p$
(Theorem~\ref{not}), being necessarily of the type $\j_1.\j_2\dots\j_r$
(Theorem~\ref{2003}). Then, using the coordinate functions of this chosen
EKR, the local equations of $S$ around $p$ (which now becomes 0) will be:
\begin{itemize}\itemsep=0pt
\item
$x_k = 0$  for all $k$ such that $j_k = 2$,
\item
$x_s = y_s = 0$  for all $s$ such that $j_s = 3$.
\end{itemize}
Indeed, Proposition \ref{posi} holds at all points. Hence, on analyzing
the key congruence (\ref{sanbis}), it is the vanishing of $x_l$ that
is decisive for the inclusion to hold true in the $l$-th sandwich.
This explains all the $x$-equations, coming from both types of letters:
the $j_k = 2$ and the $j_s = 3$\ in the code of~$\C$. In other words,
the $x$-equations are the equations of the locus of the sandwich class,
say~$\E$, which encompasses~$\C$. Now some auxiliary equations,
excising from~$\E$ the singularity class~$\C$, should be added
to them.

Continuing then, any letter $j_s = 3$ in $\C$, upon analyzing
the construction of the small f\/lag of~$D^s$, brings in the additional
equation $y_s = 0$. Such is the {\it eventual} conclusion drawn from
the arguments used in the proof of Theorem~\ref{2003}. We mean by this
that not only the germ at 0 (previously our point $p$), but precisely
the germs at all points $q$ having the coordinates~$x_\nu$,~$x_s$,~$y_s$
vanishing ($s = \nu + l + 1$ in the notation from Section~\ref{pf}),
satisfy the inclusion $V_{2l+3}(q) \subset L(D^{\nu-2})(q)$, or
$\subset F(q)$ when $\nu = 2$. That is, have a letter 3 at
the $s$-th position in {\it their} codes, preceded by $l$ letters 1,
preceded in turn by a letter not 1 at the $\nu$-th position.
(Note that, naturally, each such variable $x_\nu$ is in the union
of {\it all} variables $x_k$ and $x_s$ appearing in the present
discussion.) This explains the $y$-equations joining the previous
$x$-equations in our local description of~$S$.
\end{proof}

\section{Concrete example of discernment\\ inside the sandwich
class 1.\underline{2}.1.\underline{2}}\label{concr}

On the manifold $\R^{(4 + 1)2 + 1}(t,x_0,y_0,x_1,y_1,\dots,
x_4,y_4)$ we will propose two non-equivalent families of EKRs, both
sitting in the sandwich class 1.\underline{2}.1.\underline{2}. One is
called $D$ and has (all members of the family) the singularity class
$\W(D) = 1.2.1.2$, and the other is called $E$ and has the class $\W(E)
= 1.2.1.3$. Moreover, we will see the geometrical distinction between
$D$ and $E$ at work. This time, for bigger transparency, we use the
vertical writing of the most involved vector f\/ield's generators.

The f\/irst family of germs at $0 \in \R^{11}$ is generated
by the following vector f\/ields:
\[
D  = \Biggl(
\begin{array}{@{}r}x_4\!\!\left(\!\!\!\!\begin{array}{r}\left.x_2\!\!\left(\!\!\begin{array}{c@{}}1\\x_1\\y_1
                                                        \end{array}\right.\!\right]\\
                 \left.\begin{array}{c@{}}1\\y_2\end{array}\right]\mbox{\hskip.04mm}\\
                 \left.\begin{array}{c@{}}X_3\\Y_3\end{array}\right]\mbox{\hskip.04mm}
                 \end{array}\right.\\
\left.\begin{array}{c@{}}1\\Y_4\end{array}\right]\mbox{\hskip2mm}\\
\left.\begin{array}{c@{}}0\mbox{\hskip1mm}\\0\mbox{\hskip1mm}\end{array}\right]\mbox{\hskip2mm}
\end{array}\!\!\!\!\!\!\!,  \p_{x_4} ,  \p_{y_4}\!\!
 \Biggr).
\]
Here, a variable $x$ to the left of a bracket means (also in the sequel)
multiplying by $x$ all entries subsumed by that bracket. Thus, this is
a 3-parameter family of rank-3 distributions (some of them might be pairwise equivalent because
EKRs are only {\it pseudo}-normal forms).
And the second 2-parameter family of germs (some of them might be
pairwise equivalent as well) reads as follows:
\[
E  = \Biggl(
\begin{array}{@{}r}x_4\!\!\left(\!\!\!\!\begin{array}{r}\left.x_2\!\!\left(\!\!\begin{array}{c@{}}1\\x_1\\y_1
                                                        \end{array}\right.\!\right]\\
                 \left.\begin{array}{c@{}}1\\y_2\end{array}\right]\mbox{\hskip.01mm}\\
                 \left.\begin{array}{c@{}}X_3\\Y_3\end{array}\right]\mbox{\hskip.01mm}
                 \end{array}\right.\\
\left.\begin{array}{c@{}}y_4\\1\end{array}\right]\mbox{\hskip2mm}\\
\left.\begin{array}{c@{}}0\mbox{\hskip.9mm}\\0\mbox{\hskip.9mm}\end{array}\right]\mbox{\hskip2mm}
\end{array}\!\!\!\!\!\!\! ,  \p_{x_4} ,  \p_{y_4}
\!\!\Biggr) .
\]

{\it Attention}. These objects come directly from Theorem~\ref{2003}:
on the level of local normal forms, precisely the EKR families
\1.\2.\1.\2 and \1.\2.\1.\3 represent the sandwich
class 1.\underline{2}.1.\underline{2} which is the union of the
singularity classes 1.2.1.2 and 1.2.1.3.
\vskip1.5mm
It is straightforward (and based only on sandwich-like inclusions, cf.\
Corollary~\ref{relation}) that, in the process of constructing $\W(D)$
and $\W(E)$, both germs happen to belong to 1.\underline{2}.1.\underline{2}.
Then, in the passing from sandwich to singularity class(es), the f\/irst
\underline{2} from the left causes no trouble (see Section~\ref{refin})
while the specif\/ication of the second \underline{2} is subtler.

For that second \underline{2} in 1.\underline{2}.1.\underline{2}, the values
of the integers $\nu$ and $l$ are $\nu = 2$ and $l = 1$. Because $\nu$
takes the smallest possible value, the covariant subdistributions of $D^1$
and $E^1$ (after the standard indexation of the big f\/lags of $D$ and $E$)
enter into play. And, because we work with EKRs, these covariant objects
are both equal to
%\begin{equation}%\label{cov}
$F = \bigl(dt,dx_0,dy_0\bigr)^\perp$
%\end{equation}
(by Proposition~\ref{posi}, $F$ is spanned by all the versors
save $\p_t$, $\p_{x_0}$, $\p_{y_0}$).

 Let, for the sake of brevity, $D = D_1 \subset D_2 \subset D_3
\subset \cdots$  and  $E = E_1 \subset E_2 \subset E_3 \subset \cdots$
 be the respective small f\/lags. Since $2l + 3 = 5$, the algorithm of
f\/inding the singularity class requires to analyze the positions at 0 of
$D_5$, $E_5$, and $F$. In order to do that it is helpful to watch carefully
the early members~$D_2$ and~$E_2$:
\begin{gather*}
D_2\!  =\!  \Biggl(
\begin{array}{@{}r@{}}x_4\!\! \left(\!\! \begin{array}{@{}r}\left.x_2\!\!\left(\!\!\begin{array}{c@{}}1\\x_1\\y_1
                                                        \end{array}\right.\! \right]\\
                 \left.\begin{array}{c@{}}1\\y_2\end{array}\right]\mbox{\hskip.04mm}\\
                 \left.\begin{array}{c@{}}X_3\\Y_3\end{array}\right]\mbox{\hskip.04mm}
                 \end{array}\right.\\
\left.\begin{array}{c@{}}1\\Y_4\end{array}\right]\mbox{\hskip2mm}\\
\left.\begin{array}{c@{}}0\mbox{\hskip1mm}\\0\mbox{\hskip1mm}\end{array}\right]\mbox{\hskip2mm}
\end{array}\!\!\!\!  ,\!\!\!\! \!
\begin{array}{r}\left.x_2\!\!\left(\!\!\begin{array}{c@{}}1\\x_1\\y_1\end{array}\!\right.\right]\\
      \left.\begin{array}{c@{}}1\\y_2\end{array}\right]\mbox{\hskip.01mm}\\
      \left.\begin{array}{c@{}}X_3\\Y_3\end{array}\right]\mbox{\hskip.01mm}\\
      \left.\begin{array}{c@{}}0\mbox{\hskip1mm}\\0\mbox{\hskip1mm}\end{array}\right]
      \mbox{\hskip.3mm}\\
      \left.\begin{array}{c@{}}0\mbox{\hskip1mm}\\0\mbox{\hskip1mm}\end{array}\right]
      \mbox{\hskip.3mm}
\end{array}\!\!\!\! , \p_{x_3},    \p_{y_3}  ,  \p_{x_4}  ,  \p_{y_4} \!\!\Biggr) \!,
\ \
E_2\! = \!\Biggl(
\begin{array}{@{}r}x_4\!\!\left(\!\! \begin{array}{@{}r}\left.x_2\!\!\left(\!\!\begin{array}{c@{}}1\\x_1\\y_1
                                                        \end{array}\right.\!\right]\\
                 \left.\begin{array}{c@{}}1\\y_2\end{array}\right]\mbox{\hskip.01mm}\\
                 \left.\begin{array}{c@{}}X_3\\Y_3\end{array}\right]\mbox{\hskip.01mm}
                 \end{array}\right.\\
\left.\begin{array}{c@{}}y_4\\1\end{array}\right]\mbox{\hskip2mm}\\
\left.\begin{array}{c@{}}0\mbox{\hskip.9mm}\\0\mbox{\hskip.9mm}\end{array}\right]\mbox{\hskip2mm}
\end{array}\!\!\!\!\!\! ,\!\!\!\!\!
\begin{array}{r}\left.x_2\!\!\left(\!\!\begin{array}{c@{}}1\\x_1\\y_1\end{array}\right.\!\right]\\
      \left.\begin{array}{c@{}}1\\y_2\end{array}\right]\mbox{\hskip.01mm}\\
      \left.\begin{array}{c@{}}X_3\\Y_3\end{array}\right]\mbox{\hskip.01mm}\\
      \left.\begin{array}{c@{}}0\mbox{\hskip1mm}\\0\mbox{\hskip1mm}\end{array}\right]
      \mbox{\hskip.2mm}\\
      \left.\begin{array}{c@{}}0\mbox{\hskip1mm}\\0\mbox{\hskip1mm}\end{array}\right]
      \mbox{\hskip.2mm}
\end{array}\!\!\!\!  , \p_{x_3}  , \p_{y_3}  ,  \p_{x_4}  ,  \p_{y_4}\!\! \Biggr)\! .\!
\end{gather*}
Note that the f\/irst generators on the left in these descriptions are
superf\/luous for $D_2$ ($E_2$) as such. Yet we are to compute parts of
the small f\/lags of the departure objects ($D$ and $E$), and the presence of these generators of $D$ and $E$
makes the necessary computations easier.

 The main observation is that by multiplying the f\/irst and second
generators in $D_2$ ($E_2$)  we get  $\p_{x_2} + Y_4\p_{y_2} \in D_3$
$\bigl(y_4\p_{x_2} + \p_{y_2} \in E_3\bigr)$, compare to~(\ref{y-s}).
Because of the special role of the variable $x_2$ in $D$, this leads
in two more steps to
\[
\p_t + x_1\p_{x_0} + y_1\p_{y_0} + Y_4\p_{y_1}  \in  D_5 .
\]
This is just an instance of the inclusion (\ref{shy}) with
$Z = \p_t + x_1\p_{x_0} + y_1\p_{y_0}$, $\nu = 2$, $s = 4$, $V_5 = D_5$.
Thus $D_5(0) \not\subset F(0)$. To distinguish between
the two cases, the same operations performed for~$E$ lead to a benign
inclusion $y_4Z + \p_{y_1} \in E_5$. Also, clearly, $x_4\p_{x_2},
x_4\p_{y_2} \in E_3$, which leads to another benign inclusion $x_4Z
\in E_5$. In consequence, the inclusion $E_5(0) \subset F(0)$ holds.
(These statements are a simple instance of the main line of computations
in the proof of Theorem~\ref{2003}.) Therefore, by the def\/inition of
singularity classes, the second \underline{2} in the sandwich word
1.\underline{2}.1.\underline{2} is being specif\/ied: to~2 for~$D$,
and to 3 for $E$. And the f\/irst \underline{2} is univocally replaced
by 2.  That is to say, the 3-parameter family of germs $D$ is included
in the singularity class 1.2.1.2, whereas the 2-parameter family $E$ is
included in the class 1.2.1.3.

\section{Proof of Proposition \ref{uni}}\label{APP}
Let ${\rm rk}\,D = n + 1$  and  ${\rm cork}\,D = m > 1$.
\begin{proof} The proof is local around any f\/ixed point $p \in N$.  In
view of the Frobenius theorem there clearly exist local coordinates
$x_0, x_1,\dots, x_m, y_1, y_2,\dots, y_n$ vanishing at $p$
and such that
\begin{equation}\label{ABC}
E = \bigl(dx_0,dx_1,\dots,dx_m\bigr)^\perp=
\bigl(\p_{y_1},\p_{y_2},\dots,\p_{y_n}\bigr)
\end{equation}
and
\[
D = \bigl(\p_{x_0} + f_1\p_{x_1} + \cdots + f_m\p_{x_m},E\bigr)
\]
for certain functions $f_i$ which may even be assumed vanishing at 0:
$f_1(0) = \cdots = f_m(0) = 0$. Now the `two-step' assumption
$[D,D] = TN$ implies that $m \le n$ and
\[
{\rm rk}\,\left(\frac{\p f_i}{\p y_j}\right)(0)  =  m ,
\]
with the indices' ranges $i = 1,\dots, m$, $j = 1,\dots,n$.
Without loss of generality one can have the second range reduced
to $j = 1,\dots, m$. After such a simplif\/ication the functions
\[
\bigl(x_0,x_1,\dots,x_m,f_1,\dots,f_m,y_{m+1},\dots,
y_n\bigr)
\]
are independent at 0 and we take them as new variables (for simplicty,
we keep writing the letters $x$ and $y$ for the new variables). The
purpose is twofold. Firstly, this coordinate change mapping, say
$\phi$, is clearly of the form $\phi(x_0,x_1,\dots,x_m,\dots)
= (x_0,x_1,\dots,x_m, \dots)$, implying that the description
(\ref{ABC}) of $E$ holds in both old and new variables. Secondly,
the distribution~$D$ now gets an extremely simple description
\begin{equation}\label{BCD}
D = \bigl(\p_{x_0} + y_1\p_{x_1} + \cdots + y_m\p_{x_m},
\p_{y_1},\p_{y_2},\dots,\p_{y_n}\bigr).
\end{equation}
This in dual terms says
\[
D^\perp = \bigl(dx_1 - y_1dx_0, dx_2 - y_2dx_0,\dots,
dx_m - y_mdx_0\bigr)
\]
and allows one to easily search for the covariant object. In fact, at
each point close to $0 \in \R^{m+n+1}$ one is looking for all 1-forms
$\alpha$ such that
\begin{equation}\label{DEF}
\left.\bigl(\alpha\wedge dx_0\wedge dy_i\bigr)\right|_{D} = 0,
\qquad i = 1,2,\dots,m.
\end{equation}
In the new coordinates the answer does not depend on a point.
Indeed, upon writing
\[
\alpha = a_0dx_0 + a_1dx_1 + \cdots + a_mdx_m + b_1dy_1 + \cdots
+ b_mdy_m + \sum_{j = m+1}^n b_jdy_j,
\]
one instantly observes that the coef\/f\/icients $a_0,a_1,\dots,a_m$
are free, subject to no restrictions (for~$\left.dx_i\right|_D$,
$i = 1,\dots,m$, are multiples of $\left.dx_0\right|_D$). Concerning
the coef\/f\/icients $b_j$ with $j > m$, they vanish identically, because
the dif\/ferentials $dx_0,dy_1,\dots,dy_n$ are free also after their
restricting to $D$ (cf.~(\ref{BCD})).
As it could be expected, the key coef\/f\/icients are $b_1,\dots,b_m$.
Because~$m$ is greater than~1, the conditions~(\ref{DEF}) imply that
$0 = b_1 = \cdots = b_m$ identically\footnote{It would be otherwise
for $m = 1$, cf.~\cite[p.~165]{Mormul2004a}. The dichotomy in the
classical Lie--B\"acklund theorem has its roots precisely at this place.}.
In fact, taking $i = 1$ in~(\ref{DEF}) implies $b_2 = b_3 = \cdots =
b_m = 0$ at the point under consideration. Taking $i = 2$ implies
$b_1 = b_3 = \cdots = b_m = 0$ and already at this moment all
coef\/f\/icients $b_1,\dots,b_m$ are zero.

 That is, $\alpha = a_0dx_0 + a_1dx_1 + \cdots + a_mdx_m$ and
$a_0,a_1,\dots, a_m$ are free and this holds at every point. Therefore
at every point the $\alpha$'s describe $\bigl(dx_0, dx_1, \dots,dx_m\bigr)
^\perp = E$, and this is the covariant subdistribution $\widehat{D}$
of $D$.

Concerning the Cauchy-characteristic module $L(D)$, it is a classical
fact going back to~\cite{vonWeber} (see Theorem~1 on p.~211 there)
that for $D$ under form (\ref{BCD}) the Cauchy module is regular and
\[
L(D)  =  \bigl(dx_0, dx_1, \dots, dx_m, dy_1, \dots, dy_m\bigr)^
\perp  =  \left(\p_{y_{m+1}},\dots, \p_{y_n}\right) .
\]
Indeed, then, in view of~(\ref{ABC}), $L(D)$ is of corank~$m$
inside $E$. Proposition~\ref{uni} is now fully proved.
\end{proof}

\subsection*{Acknowledgment}
The author was supported by Polish Grant
MNSzW N\,\,N\,201\,\,397\,937.

\pdfbookmark[1]{References}{ref}
\LastPageEnding

\end{document}